\documentclass[12pt]{amsart}

\usepackage{amsmath, amssymb}
\usepackage{enumitem}
\newtheorem{theorem}{Theorem}[section]
\newtheorem{lemma}[theorem]{Lemma}
\newtheorem{prop}[theorem]{Proposition}

\theoremstyle{definition}

\theoremstyle{remark}

\numberwithin{equation}{section}

\newcommand{\rr}{{\mathbb R}}
\newcommand{\rd}{{\mathbb R^d}}
\newcommand{\zd}{{\mathbb Z^d}}

\newcommand{\N}{{\mathbb N}}
\newcommand{\nat}{{\mathbb N}}

\allowdisplaybreaks

\newcommand{\RCM}{\operatorname{RCM}}
\newcommand{\Z}{\mathbb{Z}}
\renewcommand{\P}{\mathbb{P}}
\newcommand{\Pcal}{\mathcal{P}}
\newcommand{\E}{\mathbb{E}}

\newcommand{\1}{\mathbf{1}}
\newcommand{\M}{\operatorname{M}}
\newcommand{\G}{\operatorname{G}}

\usepackage{xcolor}

\begin{document}
\sloppy
\title[Random walk on the random connection model]{Random walk on the random connection model} 

\author{Ercan S\"onmez}
\address{Ercan S\"onmez, Department of Statistics, University of Klagenfurt, Universit\"atsstraße 65--67, 9020 Klagenfurt, Austria}
\email{ercan.soenmez\@@{}aau.at} 

\author{Arnaud Rousselle}
\address{Arnaud Rousselle, Institut de Mathématiques de Bourgogne, UMR 5584 CNRS, Université Bourgogne Franche-Comté, F-2100 Dijon, France}
\email{arnaud.rousselle@u-bourgogne.fr}

\date{\today}

\begin{abstract}
We study the behavior of the random walk in a continuum independent long-range percolation model, in which two given vertices $x$ and $y$ are connected with probability that asymptotically behaves like $|x-y|^{-\alpha}$ with $\alpha>d$, where $d$ denotes the dimension of the underlying Euclidean space. More precisely, focus is on the random connection model in which the vertex set is given by the realization of a homogeneous Poisson point process. We show that this random graph exhibits the same properties as classical discrete long-range percolation models studied in \cite{B} with regard to recurrence and transience of the random walk. {The recurrence results are valid for every intensity of the Poisson point process while the transience results hold for large enough intensity.} Moreover, we address a question which is related to a conjecture in \cite{HHJ} for this graph.
\end{abstract}

\keywords{Percolation, random graphs, long-range percolation, random connection model, Poisson process, random walk in random environment, recurrence, transience}
\subjclass[2020]{Primary 05C81; Secondary 05C80, 60K35, 60G55, 82B20.}
\maketitle

\baselineskip=18pt

\section{Introduction}
The theory of random graphs is a large branch of classical as well as modern probability theory. For many decades it has been studied to a considerable extent, not least because of its increasingly important role in science. There is a variety of applications of random graph models, in particular motivated by the study of real-world networks. The reader is advised to consult the references \cite{Boll, J, VH} for an overview.

In the present paper there will be a particular interest in the study of the \textit{random connection model}. In such a graph the vertex set is obtained as a random subset of $\rd$ and is given by the realization of a homogeneous Poisson point process $\mathcal{P}$ on $\rd$ with $d \geq 1$. Given a realization of $\mathcal{P}$ an edge between two vertices $x,y \in \mathcal{P}$ is drawn with probability $g(x,y)$ depending only on the distance of $x$ and $y$ in $\rd$. Motivated by applications in communication networks, see \cite{Gil, Kes}, the random connection model has been introduced in \cite{P}, in which results about the percolation behavior and property of connected component sizes have been established. Moreover, it has further been studied in \cite{Bur, Iy, P2}, also in the case in which the underlying Poisson point process is replaced by more general stationary point processes. {We note that our main results can also be extended to more general stationary point processes provided they have a finite range of dependance and a control over the number of points that fall in boxes with high probability\,similarily as in \cite{R}.} In \cite{P} it is shown that for dimensions $d \geq 2$ percolation occurs if the intensity of the Poisson point process is sufficiently large.

The random connection model can be seen as a continuous analogue of discrete long-range percolation on the integer lattice $\zd$ introduced in \cite{Schul}, see \cite{DHH, DW} for extensions of the model and \cite{LNS} for further investigations of the random connection model. Recall the classical model of long-range percolation on $\zd$ in which two arbitrary vertices can be connected by a bond. More precisely, for every $u$ and $v$ in $\zd$ there is an edge connecting $u$ and $v$ with some probability $p(u,v)$ only depending on the distance of $u-v$ and the origin.

In this paper we focus on the behavior of a particle which performs a random walk on the graph. It is commonly used as a way of characterizing the geometry of a graph.

\textit{Random walks in random environment} have been studied for many decades and have become a fruitful research area. Its beginning is marked by the study of limit theorems for the random walk in a quite general setting, see \cite{Sol} and also \cite{Angel, Ben, GKZ, Hag}. In \cite{GKZ} the authors showed that the infinite cluster of supercritical percolation in $\zd$ is transient almost surely for all dimensions $d \geq 3$, establishing an analogy to P\'olya's theorem for the lattice $\zd$ with nearest neighbor edges. Berger \cite{B} investigated this problem in the long-range percolation model in $\zd$ in which the edge probability $p(u,v)$ between two vertices $u$ and $v$ asymptotically behaves like $|u-v|^{-\alpha}$ for some $\alpha >d$. Remarkably, he showed that the behavior of the random walk not only depends on the dimension $d$ but on $\alpha$ as well. More precisely, he showed that supercritical long-range percolation is recurrent for $d=1,2$ and $\alpha \geq 2d$, whereas it is transient in all dimensions $d \geq 1$ given that $\alpha \in (d,2d)$, which is a dramatic difference to P\'olya's theorem and the results in \cite{GKZ}. The case $\alpha \geq 2d$ in dimension $d \geq 3$ has been open ever since and it is conjectured that supercritical long-range percolation is transient in this case as well. In \cite{HHJ} the authors mention that it would be interesting to verify transience in this case.

Our main interest will be in the study of the random connection model in which we choose the edge probability $g(x,y)$ to be a function that asymptotically behaves like $|x-y|^{-\alpha}$ for $\alpha >d$. We will answer the question of recurrence and transience of the random walk in all dimensions and show that this random graph behaves similarly to long-range percolation described above. {That is we show that random walks on the random connection model are a.s.\,recurrent for $d = 2$ if $\alpha \geq 4$ (see Theorem~\ref{th:recurrent}). Furthermore, by using a suitable stochastic domination argument, we will show that, provided that the intensity of the underlying Poisson point process is sufficiently large, random walks on the infinite cluster of the random connection model are a.s.\, transient for $d=1,2$ if $\alpha \in (d,2d)$ and for $d \geq 3$ if $\alpha \geq d$  (see Theorem~\ref{transience}). Some proofs are inspired by the methods presented in \cite{B}. We note that an adaptation of the proof of \cite[Theorem 3.1]{B}, not detailed in the present paper, provides a proof of a.s.\,recurrence of the random walk on the random connnection model when $d = 1$ if $\alpha \geq 2$ from the classical Nash-William criterion (see {\it e.g.} \cite[Section 2.5]{LP}).

We also note that, in \cite{CFG}, Caputo, Faggionato and Gaudilli\`ere studied recurrence and transience properties of random walks on complete graphs generated by point processes in $\rr^d$ with jump probabilities which are a decreasing function of the distance between points. They typically choose jump probabilities that decay as $\varphi_{p,\alpha}(t)=1\wedge t^{-d-\alpha}$, $\alpha>0$ or $\varphi_{e,\beta}(t)=\exp(-t^\beta)$, $\beta>0$, where $t$ stands for the inter-points distances. Under assumptions on the point process that are in particular satisfied by the Poisson point process $\Pcal$, they show that the random walk on  $(\Pcal,\varphi_{e,\beta})$ is a.s.\, recurrent if $d=1,2$ for every $\beta >0$ and a.s.\,transient if $d\geq 3$ and $0<\beta<d$ ; they also show that $(\Pcal,\varphi_{p,\alpha})$ has a.s.\,the same type as $(\Z^d,\varphi_{p,\alpha})$, namely is transient if and only if $d=1,2$ and $0<\alpha<d$ or $d\geq 3$.} We remark that the question of transience and recurrence of random walks has also been studied for other continuum random graph models in which the underlying vertex set is given by general stationary point processes in \cite{R}. It is shown there that all the continuum random graph models in \cite{R} behave like classical nearest neighbor bond percolation as in \cite{GKZ}. However, in contrary to the present paper, the edges of all of these graphs are completely determined by the position of the vertices. {Finally, we remark that Gracar, Heydenreich, M\"onch and M\"orters announced in \cite{GHMM} recurrence/transience results for weight-dependent random connection models with scale-free degree distributions some time after the first version of this paper was submitted. 

The rest of the paper is organized as follows. In Section~\ref{sec:Def}, we introduce the random connection model under investigation and recall well known facts about random walks and electrical networks. Sections~\ref{sec:rec} and~\ref{sec:tr} are respectively devoted to the recurrence and transience results and their proofs.}

\section{Definitions and notation}\label{sec:Def}

In the sequel, we consider the Euclidean space $\rd$, $d \geq 1$, equipped with $1$-norm  $|\cdot| = |\cdot|_1$. Throughout this paper let $\mathcal{P}$ be a homogeneous Poisson point process with intensity $\rho >0$ {in $\rr^d$. This random set} is characterized by the following two properties, see {\it e.g.} \cite{Kall}:
\begin{itemize}
	\item[(i)] For every bounded set $B \in \mathcal{B} (\rd)$ the random variable $\#(\mathcal{P} \cap B)$ has a Poisson distribution with parameter $\rho \lambda_d(B)$ with $\lambda_d(B)$ denoting the Lebesgue measure of the set $B$.
	\item[(ii)] For every $n \in \nat$ and disjoint sets $B_1, \ldots, B_n \in \mathcal{B} (\rd)$ the random variables $\#(\mathcal{P} \cap B_1), \ldots, \#(\mathcal{P} \cap B_n)$ are independent.
\end{itemize}
{In the sequel, we denote by $\P_\rho$ the law of $\Pcal$ and by $\E_\rho$ the expectation w.r.t.\,$\P_\rho$.} Let also $g : \rd \to [0,1]$ be a measurable function satisfying $g(x) = g(-x)$ for all $x \in \rd$ and
\begin{equation} \label{eq:int}
 0 < \int_\rd g(x) dx < \infty.
\end{equation}

{\noindent{\bf The random connection model.} Let us recall the construction of the random connection model from \cite{P}.} 
One can interpret $\mathcal{P}$ as a process placing different particles $(X_n)_{n \in \nat}$ in $\rd$. Given such a realization for each pair $(X_i, X_j)$ of particles placed by $\mathcal{P}$ we construct a bond between $X_i$ and $X_j$ with probability $g(X_i -X_j)$, independently of all other pairs of points in $\mathcal{P}$. More formally, let $(E_{x,y}: x,y \in \rd, x \neq y)$ be a family of Bernoulli random variables. By the Kolmogorov consistency theorem{, on a suitable probability space,} we can choose $(E_{x,y}: x,y \in \rd, x \neq y)$ such that $P (E_{x,y} = 1) = g(x-y)$ for all $x,y \in \rd$ with $x \neq y$ independently{ from each other and from $\mathcal{P}$}. Given a realization of $\mathcal{P}$ and $(E_{x,y}: x,y \in \rd, x \neq y)$ we obtain the random connection model as the undirected graph with vertices given by the points $(X_n)_{n \in \nat}$ of $\mathcal{P}$ and by including an edge $(X_i, X_j)$ if and only if $E_{X_i, X_j} = 1$. We denote the corresponding graph by $\RCM (\mathcal{P})$. We denote the joint probability measure of the point process $\mathcal{P}$ with intensity $\rho >0$ and edge occupation by $P_\rho$ {and we write $E_\rho$ for the expectation under $P_\rho$}. The connected components of the resulting random graph will be called clusters. {Let us recall that} one can consider the point process $\mathcal{P}$ 'conditioned to have a point at 0' in the sense of Palm measures, see \cite{Kall}. {For Poisson point processe this is equivalent to adding} the point $X_0 =0$ to the sequence of $(X_n)_{n \in \nat}$ of $\mathcal{P}$ and as before we form a bond between $X_i$ and $X_j$, $0 \leq i < j < \infty$, with probability $g(X_i-X_j)$, independently of all other pairs $(X_i, X_j)$. {We denote by $\Pcal_0$ the Palm version of the point process and by $\RCM (\Pcal_0)$ a realization of the random connection model based on $\Pcal_0$}. Moreover, in the following we denote by $C(0)$ the cluster containing the origin. In \cite[Theorem 1]{P} it is shown that if $d \geq 2$ and \eqref{eq:int} hold there exists a critical intensity of the Poisson point process denoted by $\rho_c$ with $0 < \rho_c < \infty$ such that
$$P_\rho \big( |C(0)| = \infty \big) = 1$$
for all $\rho > \rho_c$. Throughout we will assume that $\rho > \rho_c$ always holds in dimension $d \geq 2$. We remark that in dimension $d=1$ it can be proven that the existence of an infinite cluster does not depend solely on the intensity $\rho$ anymore. Nevertheless our transience results include the case $d=1$, in particular proving the existence of an infinite cluster for this case and a particular class of functions $g$. We are ready to state a first result {that roughly means that a \emph{typical point} lies in an infinite cluster with high probability when $\rho$ is large enough.}

\begin{prop}{Assume that $d \geq 2$ and $g$ only depends on $|x|$.}
	Let $C$ be the cluster of an arbitrary vertex in $\RCM (\mathcal{P}$). Then
	$$ P_\rho \big( \#(C) = \infty \big) \to 1,$$
	as $\rho \to \infty$.
\end{prop}

\begin{proof}
	Let $M >0$ be a constant and consider a slight modification of the random connection model in which the function $g$ is replaced by the function $f$, which is given by
	$ f(x) = {g(x)} \mathbf{1}_{\{ |x| \leq M\}} $
	for every $x \in \rd$. Note that by definition the random connection model with function $g$ stochastically dominates the one with function $f$. Thus it suffices to prove the assertion in the new model. Since $f$ has bounded support, only depends on $|x|$ and is a decreasing function of $|x|$, we can apply Theorem 3 in \cite{P} and the following Corollary to obtain that
	$$ P_\rho ( \#(C) < \infty) \to 0$$
	as $\rho \to \infty$, under the new model.
\end{proof}

{A typical choice of the function $g$ is}
$$ g(x) = 1 - \exp ( -|x|^{-\alpha} ), \quad x \in \rd,$$
for $\alpha > d$. Note that the integrability condition {\eqref{eq:int}}
is satisfied, since $\alpha >d$, which is seen immediately by using the inequality $1- e^{-x} \leq x$ for small values of $x$. {We remark that an other possible choice is $$ g(x) = \mathbf{1}_{ \{|x| \leq R \} }$$
for all $x \in \rd$ and some constant $R>0$. This corresponds to the well-known \textit{Poisson blob model} (see {e.g.} \cite{Hall1, Hall2}) for which the conclusions of Theorem~\ref{th:recurrent} and Theorem~\ref{transience}~\eqref{it:tr1} hold.
\bigskip

\noindent{\bf Random walks and electrical networks}}

Let $G = (V,E)$ be a locally finite graph. Assume that each {(unoriented)} edge $\{x,y\} \in E$ is equipped with a positive conductance $c(x,y){=c(y,x)}$ and for $x \in V$ define 
$$ c(x) = \sum_{y : \{x,y\} \in E} c(x,y).$$
Recall that{, provided that $0<c(y)<+\infty$ for all $y\in V$,} a random walk on $G$ associated with the conductances of its edges is a Markov chain $(S_n)_{n \in \nat_0}$ with $S_0 \in V$ and $S_{n+1}$ is chosen at random from the neighbors of $S_n$, {\it i.e.}
$$S_{n+1} \in \big\{ x \in V: \{ x, S_n \} \in E \big\} ,$$
with probability
$$ P(S_{n+1} = x | S_n = y) = \frac{c(x,y)}{c(y)}, \quad \{x,y\} \in E,$$
independently of $S_0, \ldots, S_{n-1}$. A {connected} graph is \emph{recurrent} if for every starting point $S_0 \in V$ a random walk returns to $S_0$ almost surely { and  transient otherwise}. We note that the random connection model is a locally finite graph $P_\rho$-almost surely for all $\rho >0$. Indeed, we now show that the degree of the origin {in $\RCM (\Pcal_0)$} is finite almost surely. Let $\operatorname{deg} (0)$ be the number of vertices connected to the origin. Since the points connected to the origin can be viewed as an inhomogeneous Poisson point process we get for some constant $c$
$$\mathbb{E} [\operatorname{deg} (0) ] \leq c  \int_\rd g(x) dx < \infty$$
by the integrability condition. This implies that with probability one the point $0$ is only bonded to finitely many points and the claim follows. {It then follows from \cite[Lemma B.2]{CFP} that each vertex in $\RCM(\Pcal)$ has a.s. finite degree. In the rest of the paper, we will consider random walks on connected components of the random connection model associated with the conductances: 
\[c(x,y)=\mathbf{1}_{\{x,y\}\mbox{ is an edge of }\RCM(\Pcal)}\]
and we do not consider the degenerate cases in which a connected component is reduced to a single point. So, for almost every realization of $\RCM (\Pcal)$, the random walk is well-defined.

In the proof of Theorem~\ref{transience}, we will make use of the so-called \emph{Royden-Lyons criterion}. In order to state it, let us recall that given $x\in V$ and $B\subset V$ such that $B\not\ni x$, a \emph{unit flow from $x$ to $B$} is an antisymmetric function $f:\, V\times V\rightarrow \rr$ satisfying
\[\operatorname{div} f(y):=\sum_{z\in V}f(y,z)\left\{
\begin{array}{ll}
=1 &\mbox{if }y=x\\
=0 &\mbox{if }y\not\in B\cup\{x\}\\
\leq 0 &\mbox{if }y\in B\\
\end{array}
\right..\]
If $B=\emptyset$ then $f$ is a \emph{unit flow from $x$ to $\infty$}. The energy of $f$ is defined by
\[\mathcal{E} (f)=\frac{1}{2}\sum_{\substack{y,z\in V\\y\neq z}}\frac{f(y,z)^2}{c(y,z)}.\]
\begin{theorem}[{Royden-Lyons criterion, see \cite{L}}]\label{th:RL}

The random walk $(S_n)_{n \in \N_0}$ is transient if and only if there exists a unit flow on the electrical network from some point $x\in V$ to $\infty$ with finite energy.
\end{theorem}

For more details about random walks and electrical networks we refer {\it e.g.} to \cite{LP}.
}

\section{Recurrence in dimension 2}\label{sec:rec}
We first state the main result of this section whose proof is adapted from the one of Theorem~3.10 in \cite{B}.

\begin{theorem} \label{th:recurrent}
	Let $d = 2$. Consider the random connection model with the function $g$ decreasing with $|x|$ and satisfying
	$$ \sup_{x \in \rr^2} \frac{g(x)}{|x|^{-s}} < \infty $$
	for some $s \geq 4$.	Then, $P_\rho$-a.s., the simple random walk on any connected component of $\RCM (\mathcal{P})$ is recurrent.
\end{theorem}

In order to prove this result we will associate with almost every realization of $\RCM(\mathcal{P})$ a recurrent network on $\Z^2$ with higher effective conductance. The construction is in two steps. First, we associate with $\RCM(\mathcal{P})$ a multigraph, say $\M$, with vertex set $\Z^2$. Then, we use a \emph{projection process}, similar to the one used by Berger in \cite{B}, to obtain a new recurrent network $\G$ with only short edges in $\Z^2$, where a short edge is defined to have length of 1-norm equal to 1. Since the effective conductance increases in both steps, Theorem\,\ref{th:recurrent} follows from this construction. 

\noindent{\it Step 1.}
Let us part $\rr^2$ into boxes of side 1 : 
\[B_z=z+\left[-\frac{1}{2},\frac{1}{2}\right[^2, \qquad z\in\Z^2.\]	

Let $\RCM(\mathcal{P})$ be a realization of the random connection model. The multigraph $\M$ has vertex set $\Z^2$ and its edge set is constructed as follows. For each (unoriented) edge $(X^{(1)},X^{(2)})$ of $\RCM (\Pcal)$ with $X^{(1)}\in B_{z^{(1)}}$ and $X^{(2)}\in B_{z^{(2)}}$ for some $z^{(1)}, z^{(2)} \in \Z^2$, we add an edge of conductance 1 between $z^{(1)}$ and $z^{(2)}$. Note that $\M$ has higher effective conductance than $\RCM(\mathcal{P})$ since this procedure is equivalent to merging all the vertices in each box $B_z$, $z\in\Z$, in a single vertex at $z$ (in other words, it is equivalent to putting an infinite conductance edge between $z$ and each vertex of $\RCM (\mathcal{P})$ in $B_z$).

\noindent{\it Step 2.} We cut each edge of length $l$ in $\M$ into $l$ short edges with conductance $l$. More precisely, for $z^{(1)},z^{(2)}\in\Z^2$ with $|z^{(2)}-z^{(1)}|=l$ and $z^{(1)}_1\leq z^{(2)}_1$, let $\gamma(z^{(1)},z^{(2)})$ be the nearest neighbor path of length $l$ between $z^{(1)}$ and $z^{(2)}$ contained in $([z^{(1)}_1,z^{(2)}_1]\times \{z^{(1)}_2\})\cup (\{z^{(2)}_1\}\times [\min\{{z^{(2)}_1,z^{(2)}_2\}},\max\{{z^{(2)}_1,z^{(2)}_2\}}])$. Then, if there is an edge between $z^{(1)}$ and $z^{(2)}$ in $\M$ with $|z^{(2)}-z^{(1)}|=l$ and $z^{(1)}_1\leq z^{(2)}_1$, we erase this edge and put instead edges of conductance $l$ along $\gamma(z^{(1)},z^{(2)})$ to obtain a new multigraph $\M '$. Note that $\M '$ has higher effective conductance than $\M$. Finally, $\G$ is the network with vertices $\Z^2$ in which the conductance between $z^{(1)}$ and $z^{(2)}$ is the sum of the conductances of the edges $z^{(1)}$ and $z^{(2)}$ in $\M '$. Observe that $\G$ has the same effective conductance as $\M '$ by the parallel law.

The following lemma is in analogy with Lemma 3.8 in \cite{B}.

\begin{lemma}\label{le:rec}

Let $\G$ be the random electrical network obtained from $\RCM (\mathcal{P})$ by applying Steps 1 and 2. 

Then, if $g$ satisfies the assumption of Theorem~\ref{th:recurrent}, the following properties hold : 
\begin{enumerate}
\item \label{pt:EffCond} the effective conductance of $\,\G$ is higher or equal to the one of $\,\RCM(\mathcal{P})$;
\item \label{pt:EqDist} the conductances of edges in $\G$ are equally distributed;
\item \label{pt:finiteCond} a.s.\,only finitely many edges of $\,\RCM(\mathcal{P})$ are projected on a nearest neigbor bond of $\Z^2$; in particular, conductances in $\,\G$ are a.s.\,finite; 
\item \label{pt:L1} if $s>4$ holds, the conductance $C$ of an edge of $\G$ is $L^1$,
\item \label{pt:CauchyTail} the conductance $C$ of an edge of $\,\G$ has a Cauchy tail, {\it i.e.\,}there exists a positive constant $c$ such that $\P \left[C>cn\right]\leq n^{-1}$ for every $n \in \N$.
\end{enumerate}
\end{lemma}
\begin{proof} Claims~\eqref{pt:EffCond} and~\eqref{pt:EqDist} are satisfied by construction.

In order to prove Claim~\eqref{pt:finiteCond}, it is enough to check that the expected number of edges in $\RCM (\mathcal{P})$ that are projected on $((0,0),(0,1))$ in $\G$ is finite. Since such an edge starts in a box $B_{z^{(1)}}$ with $z^{(1)}_1\leq 0$ and ends in a box $B_{z^{(2)}}$ with $z^{(2)}_1=0$, this expected number is bounded by

\begin{align}
\sum_{\substack{z^{(1)}_1\leq 0,\, z^{(2)}_1= 0\\z^{(1)}_2\leq 0,\, z^{(2)}_2\geq 1}}
&E_\rho\left[\sum_{\substack{X^{(1)}\in\mathcal{P}\cap B(z^{(1)})\\X^{(2)}\in\mathcal{P}\cap B(z^{(2)})}}\1_{(X^{(1)},X^{(2)})\mbox{ is an egde of }\RCM(\mathcal{P})}\right]\nonumber\\&\qquad\qquad+
\sum_{\substack{z^{(1)}_1\leq 0,\, z^{(2)}_1= 0\\z^{(1)}_2\geq 1 ,\, z^{(2)}_2\leq 0}}
E_\rho\left[\sum_{\substack{X^{(1)}\in\mathcal{P}\cap B(z^{(1)})\\X^{(2)}\in\mathcal{P}\cap B(z^{(2)})}}\1_{(X^{(1)},X^{(2)})\mbox{ is an egde of }\RCM(\mathcal{P})}\right]\label{eq:ExpBound}
\end{align}
Now, observe that, for $X^{(1)}\in B(z^{(1)})$ and $X^{(2)}\in B(z^{(2)})$, $|X^{(2)}-X^{(1)}|\geq |z^{(2)}-z^{(1)}|-2$. Hence, {using the assumptions on $g$ and recalling the definitions of $P_\rho$ and $\P_\rho$ from Section~\ref{sec:Def}}, the first summand in~\eqref{eq:ExpBound} is bounded by 
\begin{align*}
&\sum_{\substack{z^{(1)}_1\leq 0,\, z^{(2)}_1= 0\\z^{(1)}_2\leq 0,\, z^{(2)}_2\geq 1}}
\E_\rho\left[\sum_{\substack{X^{(1)}\in\mathcal{P}\cap B(z^{(1)})\\X^{(2)}\in\mathcal{P}\cap B(z^{(2)})}}g\left(X^{(2)}-X^{(1)}\right)\right]\\
&\qquad\leq
\sum_{\substack{z^{(1)}_1\leq 0,\, z^{(2)}_1= 0\\z^{(1)}_2\leq 0,\, z^{(2)}_2\geq 1}}
\E_\rho\left[\sum_{\substack{x^{(1)}\in\mathcal{P}\cap B(z^{(1)})\\x^{(2)}\in\mathcal{P}\cap B(z^{(2)})}}\min\left(1,\frac{M}{(|z^{(2)}-z^{(1)}|-2)^{4}}\right)\right]\\
&\qquad=\sum_{\substack{z^{(1)}_1\leq 0,\, z^{(2)}_1= 0\\z^{(1)}_2\leq 0,\, z^{(2)}_2\geq 1}}
\min\left(1,\frac{M}{(|z^{(2)}-z^{(1)}|-2)^{4}}\right)\E_\rho\left[\#(\mathcal{P}\cap B(z^{(1)}))\#(\mathcal{P}\cap B(z^{(2)}))\right]\\
&\qquad= \sum_{\substack{z^{(1)}_1\leq 0,\, z^{(2)}_1= 0\\z^{(1)}_2\leq 0,\, z^{(2)}_2\geq 1}}\rho^2 \min\left(1,\frac{M}{(|z^{(2)}-z^{(1)}|-2)^{4}}\right)
\end{align*}
thus finite. The same computations show that the second summand in~\eqref{eq:ExpBound} is also finite. Thus, the number of projected edges is finite and, hence, this proves~\eqref{pt:finiteCond}.

Claims~\eqref{pt:L1} and~\eqref{pt:CauchyTail}  follow from similar computations.
\end{proof}

Theorem~\ref{th:recurrent} is a consequence of Lemma~\ref{le:rec} and \cite[Theorem 3.9]{B}. For the reader's convenience, we recall the latter result, but not its proof.

\begin{theorem}[{\cite[Theorem 3.9]{B}}]

Let $\G$ be a random electrical network on the nearest neighbor bonds of the lattice $\Z^2$, such that all of the edges have the same conductance distribution and this distribution has a Cauchy tail. Then, a.s., a random walk on $\G$ is recurrent. 
\end{theorem}

\section{Transience}\label{sec:tr}

The aim of this section is to prove that the infinite cluster of the random connection model is transient almost surely if $\alpha\in (d,2d)$ for $d=1,2$ and $d\geq 3$, under some assumptions on $g$, provided that the intensity of the point process is chosen sufficiently large. More precisely, the main result of this section is stated as follows.

\begin{theorem} \label{transience}

Consider the random connection model based on a Poisson point process of intensity $\rho$ with function $g$
decreasing with $|x|$ and satisfying
\[\lim_{|x|\to 0}g(x)=1.\]

\begin{enumerate}
\item\label{it:tr1} If $d\geq 3$ then there exists $\rho' \in (0, \infty)$ such that $P_\rho$-almost surely the random walk on the infinite component of $\RCM(\Pcal)$ is transient for all $\rho \geq \rho'$.
\item\label{it:tr2} If $d=1,2$ and in addition $g$ satisfies for some $c>0$ and $\alpha\in (d,2d)$:
$$ g(x) \geq 1-\exp(-c|x|^{-\alpha}), \quad x \in \rd,$$
then there exists $\rho'\in (0, \infty)$ such that $P_\rho$-almost surely the random walk on the infinite component of $\RCM(\Pcal)$ is transient for all $\rho \geq \rho'$.
\end{enumerate}
\end{theorem}

For the proofs we cover $\Z^d$ by boxes of side $2\varepsilon$ for some $\varepsilon>0$:
\[B_z=B_z(\varepsilon)= z+\left[-\varepsilon,\varepsilon\right[^d, \qquad z\in\Z^d.\]

\noindent{\bf Proof of Theorem~\ref{transience}~\eqref{it:tr1} ($d\geq 3$).}

Let $\varepsilon$ to be chosen latter. We say that an $\varepsilon$-box $B_z$ is \emph{nice} if $\Pcal\cap B_z\neq\emptyset$ and we choose a \emph{reference vertex} $X_z$ in each nice box. Consider two neighboring nice boxes $B_{z^{(1)}}$ and $B_{z^{(2)}}$ and their reference vertices $X_{z^{(1)}}$ and $X_{z^{(2)}}$. Since vertices in $B_{z^{(1)}}\cup B_{z^{(2)}}$ are within a distance at most $2(d+1)\varepsilon$, they are connected by an edge in $\RCM(\Pcal)$ with probability at least $g\left(2(d+1)\varepsilon\right)$.

Consider now the nearest-neighbor bond percolation process $(Y_e)$ on $\Z^d$ in which a bond $e=(z^{(1)},z^{(2)})$ is \emph{open} if $B_{z^{(1)}}$ and $B_{z^{(2)}}$ are nice and their reference vertices are connected by an edge in $\RCM( \Pcal)$. Then, 
\begin{align*}
P_\rho\left(e\mbox{ is open}\right)\geq g\left(2(d+1)\varepsilon\right)\left(1-\exp\left(-(2\varepsilon)^d\rho\right)\right)^2
\end{align*}
so that we can choose $\varepsilon$ small enough and then $\rho$ large enough so that the probability $p$ that an edge $e$ is open is as close to 1 as we wish. Since the bond percolation process $(Y_e)$ has a finite range of dependence, if $p$ is large enough, it dominates a supercritical bond percolation process $(Y'_e)$ on $\	Z^d$. Hence, $\RCM(\Pcal)$, $(Y_e)$ and $(Y'_e)$ can be coupled in such a way that the infinite cluster of $(Y_e)$ is, a.s., in one-to-one correspondence with a subgraph of $\RCM(\Pcal)$.
 Since the random walk on the infinite cluster of supercritical nearest-neighbor bond percolation is transient, see \cite[Theorem 1]{GKZ}, the result follows. \strut\hfill $\Box$

\bigskip

\noindent{\bf Proof of Theorem~\ref{transience}~\eqref{it:tr2}.}

Recall the model of bond-site long-range percolation on $\mathbb{Z}^d$ in which each site is open with probability $\mu < 1$, independently of all other vertices, every pair of open sites $x$ and $y$ is connected by an open edge independently with probability
$$ 1- \exp \left( - \frac{\lambda}{|x-y|^\alpha} \right) $$
for some $\lambda > 0$. {The strategy is close in spirit to the one leading to the \emph{rough embedding method} (see \cite[Theorem 2.17]{LP} and \cite{Kanai}); we do not succeed in applying directly this method. Actually, we will couple a realization of bond-site long-range percolation on $\mathbb{Z}^d$ with a suitable subgraph of almost every realization of $\RCM(\Pcal)$ for $\rho$ large enough and then make use twice of the Royden-Lyons criterion (Theorem~\ref{th:RL}). First, for almost all realizations of the long-range percolation model on $\mathbb{Z}^d$ (with suitable parameters), there exists a unit flow with finite energy from some $z$ to $\infty$ in light of the following result taken from \cite[Lemma 2.7]{B}.

\begin{lemma} \label{Noam}
	Let $d \geq 1$, $\alpha \in (d,2d)$. Consider the bond-site long-range percolation model on $\mathbb{Z}^d$ with parameters $\lambda >0$ and $0<\mu < 1$.
	
	 Then, there exists $\mu_1 < 1$ and $\lambda_1 > 0$ such that for $\lambda \geq \lambda_1$ and $\mu \geq \mu_1$ the infinite cluster on the open sites is transient.
\end{lemma}
Then, one can exhibit a unit flow with finite energy from some vertex to $\infty$ in a subgraph of $\RCM (\Pcal)$ for almost every realization. The conclusion then follows from the Royden-Lyons criterion.

To do so, we say that a $2\varepsilon$-box $B_z$ is \emph{$(\beta,M)$-good} if \begin{enumerate}
\item $\beta\leq \#(\Pcal \cap B_z)\leq M$
\item there exists $X\in \mathcal{P}\cap B_z$ with at least $\beta-1$ neighbors in $\RCM(\Pcal)\cap B_z$.
\end{enumerate}
}

 The following result is in analogy to \cite[Lemma 5.4]{HHJ}, but is slightly stronger, since it holds for all $\alpha >d$.

{
\begin{lemma} \label{tr1}
	Let $\varepsilon>0$ and assume that $g(x)$ tends to 1 when $|x|\rightarrow 0$.  
	
	The random variables $\1_{B_z\mbox{ is $(\beta,M)$-good}}$, $z\in \Z^d$, are independent.
	Moreover, for every $\mu \in [0,1)$ and $\beta > 0$, there exists $\hat{\rho} >0$ such that for all $\rho > \hat{\rho}${, if $M=M(\rho)$ is large enough,} for any $z\in\Z^d$:
	\begin{equation}\label{eq:ProbGoodBox}
	\P\left( B_z\mbox{ is ${(\beta,M)}$-good} \right) \geq \mu .
\end{equation}	  
	
	In particular, the process $(\1_{B_z\mbox{ is $(\beta,M)$-good}})_{z\in \Z^d}$ stochastically dominates an independent site percolation process in $\Z^d$ with parameter $\mu$ provided that $\rho$ {and $M$ are} large enough. 
\end{lemma}

\begin{proof}
	Note that, in order to decide if $B_z$ is $\beta$-good or not, it suffices to know $\Pcal\cap B_z$ and which pair of points of this box are connected by an edge. The independence of the random variables $\1_{B_z\mbox{ is ${(\beta,M)}$-good}}$, $z\in \Z^d$, is thus trivial.
	
	By translation invariance, it suffices to show that~\eqref{eq:ProbGoodBox} holds for $z=0$ and $\rho$ large enough. For $\varepsilon'\in (0,\varepsilon)$ to be chosen latter, let $F_{\varepsilon'}$ be the event that there exists $X\in \mathcal{P}\cap [-\varepsilon', \varepsilon']^d$ with at least $\beta-1$ neighbors in $\RCM(\Pcal)\cap [-\varepsilon', \varepsilon']^d$ and for $k \in \N$ let $A_k$ be the event that $[-\varepsilon, \varepsilon]^d$ contains $k$ Poisson points. Then
	
	$$ P_\rho\left( B_z\mbox{ is ${(\beta,M)}$-good} \right)\geq  \sum_{k = \beta}^M P_\rho\left(F_\varepsilon | A_k\right) \\P_\rho(A_k) $$
	and, by independence, for $k \geq \beta$ we can estimate
	\begin{align*}
	P_\rho (F_{\varepsilon'} | A_k)  \geq \left(\inf_{x\in [-\varepsilon',\varepsilon']^d}g(x)\right)^{\beta-1} 
	\end{align*}
	which is as close to $1$ as we wish by choosing $\varepsilon' >0$ sufficiently small for every given $\beta \geq 1$. By definition the number of Poisson points lying in $[-\varepsilon', \varepsilon']^d$ has a Poisson distribution with parameter $\rho (2\varepsilon')^d$. Given the above choice of $\varepsilon'$ we can choose $\rho >0$ large enough such that the probability that $[-\varepsilon', \varepsilon']^d$ contains more than $\beta$ Poisson points is as close to $1$ as we wish. {It then remains to chose $M=M(\rho)$ large enough so that $\P_\rho(\#(\P\cap B_z))$ is as close to 1 as needed.} This closes the proof.
\end{proof}

The previous Lemma will allow us to control the site percolation part in the assumptions of Lemma~\ref{Noam}, namely the process of good boxes stochastically dominates an independent site percolation process on $\Z^d$ with parameter $\mu\geq \mu_1$ if $\rho$ is large enough, for any given choice of $\varepsilon,\beta>0$. For the assumption on the probability that an edge is open we will need the following result.
 Its proof is straightfoward.

\begin{lemma} \label{tr2}
Let $\varepsilon>0$ and $z^{(1)},z^{(2)}\in\Z^d$ with $|z^{(2)}-z^{(1)}|=k>0$. Let $C_1$ (resp. $C_2$) be a cluster contained in $B_{z^{(1)}}$ (resp. $B_{z^{(2)}}$). Assume that $g$ decreases with $|x|$ and satisfies for some $c>0$ and $\alpha >d$: 
    \[g(x)\geq 1-\exp\left(-c|x|^{-\alpha}\right). \]

 Then, there exist $\kappa>0$ such that
	\begin{align*}
		P_\rho\left( C_1 \text{ is connected by an open edge to } C_2 \big| \#(C_1), \# (C_2) \geq \beta  \right)
		 \geq 1- \exp \left( -\kappa \frac{\beta^2}{k^\alpha} \right) .
	\end{align*}
\end{lemma}

Let $\mu_1$ and $\lambda_1$ as in Lemma~\ref{Noam}. For $\varepsilon>0$ fixed, we first choose $\beta=\sqrt{\lambda_1/\kappa}$ and then $\rho$ {and $M$} large enough so that the process of good boxes stochastically dominates an independent site percolation process with parameter $\mu_1$ by Lemma~\ref{tr1}. In each good box $B_z$, we select a vertex $X_z$ with at least $\beta-1$ neighbors in $\RCM(\Pcal)\cap B_z$ and call it the \emph{reference vertex} of the box. Together with its neighbors in the box, it forms the \emph{reference cluster} $C_z$ of $B_z$. Then, { consider the subgraph $G$ of $\RCM(\Pcal)$ which consists in the reference clusters $C_z$ of good boxes together with the edges with endpoints in these clusters. Note that, by Lemma \ref{tr2}, if $B_{z^{(1)}}$ and $B_{z^{(2)}}$ are two distinct good boxes, their reference clusters are more likely to be connected than $z_1$ and $z_2$ needed in Lemma \ref{Noam}. By using a similar approach as in the proof of \cite[Theorem]{B} (see also \cite{CFG}), one can construct a unit flow from some vertex in $G$ to $\infty$ from a unit flow with finite energy from some $z$ to $\infty$ in the infinite cluster of bond-site long-range percolation with parameter $\mu_1$ and $\lambda_1$. To ensure that this flow has finite energy, it remains to argue that the energy is uniformly bounded in each $(\beta, M)$-good box. This is the case because there are, by definition, at most $M$ points in a $(\beta, M)$-good box. This closes the proof.} \strut\hfill $\Box$

\section*{Acknowledgements}
The first author thanks Markus Heydenreich for presenting these problems to him and for valuable discussions. The IMB receives support from the EIPHI Graduate School (contract ANR-17-EURE-0002).

\bibliographystyle{plain}

\begin{thebibliography}{10}

\bibitem{Angel} O. Angel, I. Benjamini, N. Berger, Y. Peres, Transience of percolation clusters on wedges,  Electron. J. Probab. 11, 655--669 (2006).

\bibitem{Ben} I. Benjamini, R. Pemantle, Y. Peres, Unpredictable paths and percolation, Ann. Probab. 26, 1198--1211 (1998).

\bibitem{B} N. Berger, Transience, Recurrence and Critical Behavior for Long-Range Percolation, Commun. Math. Phys. 226, 531 -- 558 (2002). Corrected proof of Lemma 2.3 available at arXiv:math/0110296v3. MR1896880.

\bibitem{Bur} R. Burton, R. Meester, Long range percolation in stationary point processes, Random
Structures Algorithms 4, 177--190 (1993).

\bibitem{Boll} B. Bollob\'as,  Random Graphs, 2nd edition, Cambridge University Press, Cambridge (2001).

\bibitem{CFG} P. Caputo, A. Faggionato, A. Gaudilli\`ere,
              Recurrence and transience for long range reversible random
              walks on a random point process, Electron. J. Probab. 14, no. 90, 2580--2616 (2009).

\bibitem{CFP} P. Caputo, A. Faggionato, T. Prescott, Invariance principle for {M}ott variable range hopping and other walks on point processes, Ann. Inst. Henri Poincar\'{e} Probab. Stat. 49 (3), 654--697 (2013).

\bibitem{DHH} M. Deijfen, H. van der Hofstad, G. Hooghiemstra, Scale-free percolation. Ann. Inst.
Henri Poincar\'e Probab. Stat. 49, 817--838 (2013).

\bibitem{DW} P. Deprez, M. V. W\"uthrich, Networks, random graphs and percolation, Theoretical
aspects of spatial-temporal modeling, 95--124, Springer, Tokyo (2015).

\bibitem{Gil} E. N. Gilbert, Random plane networks, J. Soc. Ind. Appl. Math. 9, 533--543 (1961).

\bibitem{GHMM} P. Gracar, M. Heydenreich, C. M\"onch, P. M\"orters, Recurrence versus Transience for Weight-Dependent Random Connection Models, Preprint available at https://arxiv.org/abs/1911.04350.

\bibitem{GKZ} G. Grimmett, H. Kesten, Y. Zhang, Random walk on the infinite cluster of the percolation model, Probab. Th. Rel. Fields 96, 33--44 (1993).

\bibitem{Hag} O. H\"aggström, E. Mossel, Nearest-neighbor walks with low predictability profile and percolation in $2 + \varepsilon$ dimensions, Ann. Probab. 26, 1212--1231 (1998).

\bibitem{Hall1} P. Hall, On continuum percolation, Ann. Prob. 13, 1250--1266 (1985).

\bibitem{Hall2} P. Hall, Clump counts in a mosaic, Ann. Prob. 14, 424--458 (1986).

\bibitem{HHJ} M. Heydenreich, T. Hulshof, J. Jorritsma, Structures in Supercritical Scale-free percolation, Ann. of Appl. Prob. 27, 2569--2604 (2017). 

\bibitem{Iy} S.K. Iyer, The random connection model: connectivity, edge lengths, and degree distributions,
Random Structures Algorithms 52, 283--300 (2018).

\bibitem{J} S. Janson, T. Luczak, A. Rucinski, Random Graphs, Wiley-Interscience, New York (2000).

\bibitem{Kall} O. Kallenberg, Random Measures, Akademie-Verlag, Berlin, (1986).

\bibitem{Kanai} M. Kanai, Rough isometries and the parabolicity of {R}iemannian
              manifolds, Journal of the Mathematical Society of Japan, 38 (2), 227--238 (1986).
    
\bibitem{Kes} H. Kesten, Percolation theory and first-passage percolation, Ann. Prob. 15, 1231--1271 (1987).

\bibitem{LNS} G. Last, F. Nestmann, M. Schulte, The random connection model and functions of edge-marked Poisson processes: second order properties and normal approximation, Preprint available at https://arxiv.org/abs/1808.01203.


\bibitem{LP} R. Lyons, Y. Peres, Probability on Trees and Networks,
   Cambridge Series in Statistical and Probabilistic Mathematics,
    Vol. 42, Cambridge University Press, New York (2016).
    
\bibitem{L} T. Lyons, A simple criterion for transience of a reversible {M}arkov
              chain, Ann. Probab. 11 (2), 393--402 (1983).
    
\bibitem{P} M. D. Penrose, On a continuum percolation model, Adv. Appl. Prob. 23, 536--556 (1991).

\bibitem{P2} M. D. Penrose, Connectivity of soft random geometric graphs, Ann. Appl. Probab. 26, 986--1028. (2016).

\bibitem{Peres} Y. Peres, Probability on trees: An introductory climb, In: Lectures on Probability Theory and Statistics, Lecture Notes in Math. 1717, Springer, Berlin Heidelberg, 193--280 (1993).

\bibitem{Resnick} I. Resnick, Extreme Values, Regular Variation and Point Processes, Springer-Verlag, New
York (1987).

\bibitem{R} A. Rousselle, Recurrence or transience of random walks on random graphs generated by point processes in $\rd$, Stoch. Proc. Appl. 125, 4351--4374 (2015).

\bibitem{Schul} L.S. Schulman, Long-range percolation in one dimension, J. Phys. A 16, no. 17, L639--L641 (1983).

\bibitem{Sol} F. Solomon, Random walks in random environment, Ann. Prob. 3, 1--31 (1975).

\bibitem{VH} R. van der Hofstad, Random Graphs and Complex Networks. Vol. 1, Cambridge University
Press, Cambridge (2017).


\end{thebibliography}

\end{document}